\documentclass[reqno,10pt]{amsart}
\usepackage{amsmath}
\usepackage{fullpage}
\usepackage{microtype}
\usepackage{hyperref}
\newtheorem{global-theorem}{Theorem}
\newtheorem{theorem}{Theorem}[section]
\newtheorem{lemma}[theorem]{Lemma}

\newtheorem{proposition}[theorem]{Proposition}
\newtheorem{conjecture}[theorem]{Conjecture}
\newtheorem{problem}[theorem]{Problem}
\theoremstyle{definition}
\newtheorem{definition}[theorem]{Definition}

\newtheorem{remark}[theorem]{Remark}
\usepackage[parfill]{parskip}

\begin{document}

\title{A Computational Method for Large Scale Searching of Counterexamples to the Strong Coincidence Conjecture}
\author{Scott Balchin}
\address{Department of Mathematics\\
University of Leicester\\
University Road, Leicester LE1 7RH, England, UK}
\email{slb85@le.ac.uk}

\begin{abstract}
In this short note we report on results on a computational search for a counterexample to the strong coincidence conjecture.  In particular, we discuss the method used so that further searches can be conducted.
\end{abstract}

\maketitle

\vspace{-10mm}

\section*{Introduction}

In \cite{MR3622260} a search for counterexamples for the strong coincidence conjecture on substitutions on $3$ and $4$ letters was conducted (a similar search was also carried out in \cite{MR3295586}).  The method used  was to generate substitutions based on some ordering, and then check their abelianisation for the irreducible Pisot condition.  If the substitution had the property, we then checked for strong coincidence. The main difficulty in this approach is checking that the matrix is irreducible Pisot.  Asserting that a polynomial is irreducible is a difficult problem, and also being confident that no eigenvalues lie exactly on the unit circle is non-trivial.  For 2, 3, and 4 letters we can solve these problems exactly using combinatorial methods, however for $>4$ letters we cannot do this (lack of solution by radicals etc.).

The advantage of the method that we will introduce in this report is that it is slightly easier to scale up to larger alphabets.  As a result, we were able to check 130,000 8-letter substitutions (no counterexample was found during this search).  Searches on larger alphabets are possible, however we outline some issues that need resolving for this to be accomplished.

Full commented code that was used to carry out the search is available at:

\vspace{\parskip}

\vspace{\parskip}

\centerline{\url{https://github.com/GroutGUI/PisotCheck}}

\section{Description of Approach}

We wish to search for counterexamples to the following conjecture.

\begin{conjecture}
If  a substitution is irreducible Pisot then it is strongly coincident.
\end{conjecture}

The roadmap for the method is the following:

$$\text{Find irreducible Pisot matrices using a sufficient condition}
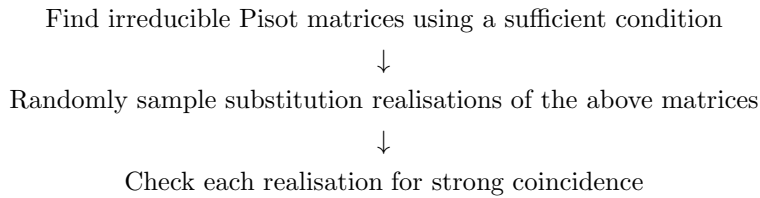
$$
$$\downarrow$$
$$\text{Randomly sample substitution realisations of the above matrices}$$
$$\downarrow$$
$$\text{Check each realisation for strong coincidence}$$

As we mentioned previously, checking if an arbitrary large matrix is irreducible Pisot is a hard problem.  To circumvent this problem we will check matrices for a sufficient (not necessary) condition.

\begin{proposition}\label{perron}
Let $A$ be a primitive matrix with characteristic polynomial $\lambda^n + a_{n-1}\lambda^{n-1} + \cdots + a_{1}\lambda + a_0$.  If:
\begin{itemize}
\item $a_0 \neq 0$,
\item $|\text{tr}(A)| = |a_{n-1}| > 1 + \sum\limits_{i = 0}^{n-2} |a_i|$,
\end{itemize}
then $A$ is irreducible Pisot.
\end{proposition}

\begin{remark}
The above condition is indeed only sufficient.  For non-necessity, consider the Tribonacci substitution (known to be irreducible Pisot) whose abelianisation has characteristic polynomial $\lambda^3-\lambda^2-\lambda$.  
\end{remark}

This condition is easy to check as we only need the coefficients of the characteristic polynomial, and  as we are working with integer matrices, using the iterative Faddeev-LeVerrier method suffices \cite{MR1213180}.  Note that in addition to this condition, we must also check that the matrix is primitive, which is not computationally intensive. In particular, by a result of Wielandt, we know a primitive $n \times n$ matrix $A$ satisfies the condition that $A^{n^2-2n + 2}$ is totally positive \cite{MR0035265}.

The first step in the overall algorithm is to find matrices which have this condition.  We need a way to generate matrices which have the conditions of  Proposition \ref{perron}.  For this experiment we chose to randomly generate binary $n \times n$ matrices where each element is idd Bernoulli with $p=0.5$.  We then check Proposition \ref{perron} and primitivity for each matrix.

In practice, this method is not efficient, but does obtain results.  In a check of $\approx$5.8bn 8x8 matrices generated in this way, only 20 satisfied the required properties.  In a check of $\approx$50bn 9x9 matrices, no such matrix was found.  Therefore, it would be useful to have a better intuition on how to generate matrices with the required constraints.

\begin{problem}
Is there a better way to randomly generate matrices to ensure a higher chance of finding irreducible Pisot ones?
\end{problem}

Given a matrix, one can ask which substitutions have that matrix as its abelianisation.  This will be important for us, as we wish to know which substitutions are irreducible Pisot.

\begin{definition}
Given a positive $n \times n$ integer matrix $A$, a \emph{realisation} of $A$ is a substitution $\phi \colon \mathcal{A} \to \mathcal{A}^+$ with $\mathcal{A} = \{a_1, \dots, a_n\}$ such that $M_\phi = A$.  We will denote the set of realisations as $\mathcal{R}(A)$.
\end{definition}

\begin{lemma}
Let $A = (a_{ij})$ be an $n \times n$ positive integer matrix.  Then
\[
|\mathcal{R}(A)| = \prod_j \dfrac{\left( \sum\limits_i a_{ij} \right)!}{\prod\limits_i \left( a_{ij} ! \right)}
\]
\end{lemma}

To put the above computation into perspective, for a random 8x8 binary matrix as generated above, $E(|\mathcal{R}(A)|) \approx 10^{13}$.  Therefore, it is not unreasonable to randomly sample a single matrix $A$ a large amount of times, as the probability of running into conjugate substitutions is minimal.

For each one of these substitutions we then check for strong coincidence.  Note that the algorithm we use for strong coincidence terminates iff the substitution is strongly coincident.  Therefore we set a limit on the number of iterations that the algorithm can take.  In practice, this limit has never been met, but if a substitution were to fail the check, then it should be rerun with a higher limit to make sure.  This problem could be eradicated if a solution to the following problem exists.

\begin{problem}
Does there exist a provable bound on the amount of iterations needed for a strong coincidence to occur?
\end{problem}

\section{Results}

1,300 8x8 irreducible Pisot matrices were found.  For each one of these matrices we sampled 100 substitutions and checked for strong coincidence.  Of these 130,000 substitutions, all passed the strong coincidence test.  Therefore no counterexample has been found.  The entire process was run on the ALICE High Performance Computing Facility at the University of Leicester, and took $\approx 24$ hours. Note that the job was split in parallel over many processors due to its stochastic nature.  In particular, the \texttt{MATLAB} code used 100 CPU hours and the \texttt{C++} code took 260 CPU hours.

\section*{Acknowledgements}

The idea of conducting a computational search for the strong coincidence conjecture in this way was kindled through many fruitful discussions with the participants of the  Workshop on Aperiodic Order at Bielefeld University (May 2017).

\bibliographystyle{siam}
\bibliography{pisot}

\end{document}